\documentclass{article}
\usepackage[top=3cm,bottom=3cm,left=3.2cm,right=3.2cm]{geometry}
\usepackage[dvipsnames,svgnames,x11names]{xcolor}
\usepackage{graphicx}
\usepackage{amsmath,amssymb,amsfonts}
\usepackage{float}
\usepackage{setspace}
\usepackage{parskip}

\begin{document}

{\Large

\thispagestyle{empty}
\vspace*{2cm}

\begin{center}
{\huge\bf{On the Structure of Asymptotic}}\\[6pt]
{\huge\bf{Space of the Lobachevsky Plane}}
\end{center}

\vspace{0.8cm}

\begin{center}
A. Shnirelman (Concordia University, Montreal, Canada)\footnote{alexander.shnirelman@concordia.ca}
\end{center}

\vspace{0.6cm}

\vspace{1.2cm}

\begin{abstract}
The notion of asymptotic space for an unbounded metric space has been introduced by Micha Gromov in 1980s. It is intended to capture the structure of a metric space at infinity. The most comprehensive definition of asymptotic space is given in the lahguage of Nonstandard Analysis (NSA). It turns out that the asymptotic space depends on the underlying nonstandard extension of the standard universe. This paper contains the exhaustive description of asymptotic spaces of the Lobachevsky plane which turns ourt to be an R-tree. However, there turn out to be a plenty of different nonisometric asymptotic spaces, including the spaces of high cardinality.
\end{abstract}

\bigskip

\section*{1 \quad Introduction}

\noindent\textbf{1.} The notion of asymptotic space, or asymptotic cone , has been introduced
by M. Gromov [5] and is intended to capture the structure of metric spaces
``at infinity''. Let $\mathcal{X} = (X,d)$ be a metric space with the set of points $X$
and the distance function $d(x,y)$. Assume that the space $\mathcal{X}$ is \textit{unbounded},
or has infinite diameter; this means that $\sup_{x,y} d(x,y) = \infty$. In the original
definition it is also required that $\mathcal{X}$ is a geodesic space, i.e. for every two
points $x,y \in X$ there is an isometric embedding $i : [0,a] \to X$, such that
$i(0) = x$, $i(a) = y$, $a = d(x,y)$; but we shall not impose this restriction,
having in mind different examples of metric spaces.

\setlength{\parindent}{1.5em}

Intuitively, asymptotic space is described as follows. Let $\mathcal{X}_{\varepsilon} = (X, d_{\varepsilon})$
be another metric space with the same set of points $X$ and a new distance
function $d_{\varepsilon}(x,y) = \varepsilon \cdot d(x,y)$, $\varepsilon > 0$. Thus, the space $\mathcal{X}_{\varepsilon}$ is obtained from
the space $\mathcal{X}$ by ``shrinking''. The asymptotic space $\mathcal{X}_0$ is defined as a ``limit'' 
 of the spaces $\mathcal{X}_{\varepsilon}$ when $\varepsilon \to 0$. Surely, we have to define what meaning the
word ``limit'' has. For example, we may define it in the sense of the Hausdorff
distance between the metric spaces. But for unbounded spaces the Hausdorff
distance is usually infinite, and this definition is restricted to a very narrow
class of metric spaces.

Other possibility to describe the asymptotic structure is the notion of
``asymptotic subcone'' [6], also introduced by M. Gromov. The space $Z =
(Z, d_Z)$ is called asymptotic subcone of the space $\mathcal{X}$, if for every finite set
of points $(z_1, \ldots, z_n) \subset Z$ there exist a sequence $\varepsilon_i \to 0$ and a sequence of
$n$-tuples $(x_1^i, \ldots, x_n^i) \subset X$, such that $\lim_{i \to \infty} \varepsilon_i d(x_j^i, x_k^i) = d_Z(z_j, z_k)$, $j, k =
1, \ldots, n$.

We see that every subset of $Z$ is also an asymptotic subcone; so, asymptotic subcones give approximation \textit{from below} to the ``true'' asymptotic space.
Every asymptotic subcone gives an idea how ``rich'' is the structure of $\mathcal{X}$ at
infinity, but they may never capture the whole structure.

\noindent\textbf{2.} The most powerful and at the same time flexible definition of asymptotic
space is formulated in the terms of nonstandard analysis. In a few paragraphs
we shall remind what the nonstandard analysis is.

Consider a set $S_0$, which will be called a \textit{basic set}; for our purposes $S_0$
may be the set $\mathbf{R}$ of real numbers. Let us form the sets $S_i$, $i = 1, 2, \cdots$,
by the inductive rule $S_{i+1} = S_i \cup \mathcal{P}(S_i)$, where $\mathcal{P}(S)$ denotes the set of all
subsets of $S$. Let $\mathcal{S} = \cup_{i \in \mathbf{N}} S_i \cup S_0$; thus, $\mathcal{S}$ is a set having some cardinality
$|\mathcal{S}|$. $\mathcal{S}$ is called a \textit{superstructure} built on the basic set $S_0$. Elements of $S_0$ are
called \textit{individuals}, and elements of $\mathcal{S} \setminus S_0$ are called \textit{entities}. It is easy to see
that all objects studied in analysis (functions, spaces, etc.) may be regarded
as elements of $S_i$ for some finite $i$, i.e. entities; thus, $\mathcal{S}$ is an appropriate
universum for all our constructions (see [4] for more details).

All first-order properties of $\mathcal{S}$ may be expressed in a language $\mathcal{L}$, whose
list of symbols consists of countable number of variables $x, y, z, \cdots$, logical
symbols $\equiv, \lor, \land, \neg, \rightarrow, \leftrightarrow$, quantors $\forall, \exists$, left and right brackets, $=$ and $\in$ as
symbols of binary relations in $\mathcal{S}$, and a set of constant symbols $c_i$, $i < \alpha$,
where the cardinal $\alpha$ is sufficiently big, $|\alpha| > |\mathcal{S}|$, so that we have enough
constants to give different names to all individuals and entities. If we assign
fixed constants as names of all elements of $\mathcal{S}$, leaving other constants idle,
we may compose all finite formulas in $\mathcal{L}$, following some simple rules [4], and
pick those of them which are true in $\mathcal{S}$ in the usual set-theoretical sense.

These formulas should be bounded; this means that all variables entering
the formula should belong to some entity. The set $T_{\mathcal{S}}$ of these formulas is a
\textit{theory} of $\mathcal{S}$, and $\mathcal{S}$ is a (standard) \textit{model} of $T_{\mathcal{S}}$.

\indent The theory $T_{\mathcal{S}}$ contains as its fragments all particular theories of analysis;
they differ by entities involved, depending on them as parameters. In particular, it contains the theory of every concrete metric space $\mathcal{X}$ (every such
space is an entity), and the theory of all metric spaces of bounded cardinality
(the set of such spaces is also an entity).

\indent Suppose now that $M = (A, E, I)$ is a set $A$ with binary relations $E$, $I$.
Suppose that $f : \mathcal{S} \to A$ is an injective map; we may regard $\mathcal{S}$ as a subset
of $A$, and $A$ as an extension of $\mathcal{S}$. Suppose that the relations $E$ and $I$ are
extensions of relations $=$ and $\in$ from the set $\mathcal{S}$ to a larger set $A$. Then every
formula in the language $\mathcal{L}$ may be interpreted in the set $A$, if we interpret
every constant, which is a name of some $x \in \mathcal{S}$, as a name of $f(x) \in A$, and
the relations $=$ and $\in$ as $E$ and $I$. If every formula from $T_{\mathcal{S}}$ is true in $A$ with
this interpretation (and every formula, which does not enter $T_{\mathcal{S}}$, is false in $A$,
for if some formula is not true, then its negation is true), then $M = (A, E, I)$
is called a \textit{model} of $\mathcal{S}$ [3]. The map $f$ is usually denoted by ${}^*$; so, $f(x) = {}^*x$
for every $x \in \mathcal{S}$. If $P(x_1, \cdots, x_n)$ is a relation in $\mathcal{S}$, i.e. an entity ($P$ is a
subset in the Cartesian product of $n$ entities; so, $P$ is an entity), then ${}^*P$ is
a relation in $A$, because the fact that the elements $y_1, \cdots, y_n \in A$ satisfy it
may be expressed by the same formula as in $\mathcal{S}$ with new interpretation. In
particular, the ${}^*$-image of a function is a function. The ${}^*$-images of entities
of $\mathcal{S}$ are called \textit{standard objects} of $A$; all other elements of $A$ are called
\textit{nonstandard}.

\indent Consider, in particular, ${}^*\mathbf{R}$, the ${}^*$-image of the basic set $S_0 = \mathbf{R}$ (which
is an entity, for $\mathbf{R} = S_0 \in S_1 \subset \mathcal{S}$). If there exist elements $y \in {}^*\mathbf{R}$, which
are not of the form ${}^*x$, $x \in \mathbf{R}$, then the model $M$ is called \textit{nonstandard};
in this case ${}^*\mathbf{R}$ is an ordered field, which is not archimedean (the axioms of
an ordered field have the form of bounded formulas of the language $\mathcal{L}$, and
therefore they are transferred to the model $M$; but the archimedean property
cannot be expressed by such a formula). Thus, ${}^*\mathbf{R}$ contains infinitely large
and infinitesimally small elements.

\indent In what follows we assume that $M$ is a nonstandard model of the universum $\mathcal{S}$.

\textbf{3.} Let us return to the definition of an asymptotic space. Let $\mathcal{X} = (X, d)$ be a metric space. Assume that $|X| < |\mathcal{S}|$; then $X$ may be embedded in $S_i$ for
some finite $i$, $d$ is a function on $X \times X$, i.e. $d \in S_{i+5}$, and $\mathcal{X} \in S_{i+7} \in \mathcal{S}$ (see
[4] for similar calculations).

\indent Let $A$ be a nonstandard model of $\mathcal{S}$. Let ${}^*\mathcal{X}$ be a ${}^*$-image of the metric
space $\mathcal{X}$ in $A$; ${}^*\mathcal{X} = ({}^*X, {}^*d)$, where ${}^*X$ may be regarded as a set in the
sense that for some elements $y$ of $A$, $yI{}^*X$; ${}^*d$ is an ${}^*\mathbf{R}$-valued function on
${}^*X \times {}^*X$. The space $\mathcal{X}$ is canonically and isometrically embedded in ${}^*\mathcal{X}$. Let
us fix a point $O \in X$ (its nonstandard counterpart ${}^*O \in {}^*\mathcal{X}$). Let us pick
some infinitesimally small $\varepsilon \in {}^*\mathbf{R}$, and define a new distance function in ${}^*\mathcal{X}$,
${}^*d_\varepsilon(x,y) = \varepsilon \cdot {}^*d(x,y)$. Thus we have obtained a new non-standard metric
space ${}^*\mathcal{X}_\varepsilon$. This space is unbounded in the following sense: for every $r \in {}^*\mathbf{R}$
there exist two points $x, y \in {}^*X$, such that ${}^*d(x,y) > r$. In fact, the original
standard space ${}^*\mathcal{X}$ is unbounded by our assumption; by the transfer principle
of the Nonstandard Analysis [4], ${}^*\mathcal{X}$ is unbounded in the nonstandard sense,
i.e. for every $R \in {}^*\mathbf{R}$ there exist points $x, y \in {}^*\mathcal{X}$, such that ${}^*d(x,y) > R$.
Taking $R = r/\varepsilon$, we find the points $x, y \in {}^*X$, such that ${}^*d_\varepsilon(x,y) > r$.

\indent Now, consider the subset $Y \subset {}^*X$, such that $y \in Y$ if and only if ${}^*d_\varepsilon(O,y)$
is finite. $Y$ is an external set; this means that $Y$ cannot be defined in the
language ${}^*L$ [4]. $Y$ may be called a finite part of ${}^*\mathcal{X}_\varepsilon$. For every two points
$x, y \in Y$ define the new distance $\delta(x,y) = \operatorname{st} {}^*d_\varepsilon(x,y)$. Here $\operatorname{st}(a)$ is a
standard part of $a$, which is defined for finite $a \in {}^*\mathbf{R}$ as a unique standard
number $a_0$, such that $a - a_0$ is infinitesimally small.

\indent Let us say that $x \approx y$, if $x, y \in Y$, and $\delta(x,y) = 0$. The factor-space
$Y/{\approx}$ is a standard metric space $X_0$, with the distance function $\delta$, which we
denote by $\mathrm{st}({}^*\mathcal{X}_\varepsilon)$. (Warning: it is not always true that ${}^*(\mathrm{st}({}^*\mathcal{X}_\varepsilon)) = {}^*\mathcal{X}_\varepsilon$;
see examples below).

\textbf{Definition 1.1} \textit{The space $\mathcal{X}_0 = \mathrm{st}({}^*\mathcal{X}_\varepsilon)$ is called the asymptotic space of the
metric space $\mathcal{X}$.}

(We always understand that the asymptotic space $\mathcal{X}_0$ depends on the
nonstandard model $M$ and the infinitesimal $\varepsilon \in {}^*\mathbf{R}$).

\indent Note that this definition is in a good correspondence with the definition
of limit in the nonstandard analysis. For example, if $f(x)$ is a real function defined for $x > 0$, then $\lim_{x \to 0} f(x)$ is defined as follows: consider the
standard universum $\mathcal{S}$, containing real numbers and functions; let ${}^*\mathcal{S}$ be its
nonstandard model, such that ${}^*\mathbf{R}$ is a proper extension of $\mathbf{R}$. Let $\varepsilon$ be an infinitely small positive nonstandard number, and $^*f(x)$ be a nonstandard counterpart of $f(x)$. This is a function from $^*\mathbf{R}$ into $^*\mathbf{R}$. In particular, $^*f(\varepsilon)$ is defined. If $^*f(\varepsilon)$ is finite, then we may define its standard part, $\mathrm{st}(^*f(\varepsilon))$. This is, \textit{by definition}, the limit of $f(x)$, when $x \to 0$.

If the limit exists in the usual sense, then we obtain the same standard number: $\lim_{x\to 0}f(x) = \mathrm{st}(^*f(\varepsilon))$. But the right hand side is defined for every bounded function, for example, for the function $f(x) = \sin 1/x$; its value depends on the non-standard model and on the particular choice of the infinitesimal $\varepsilon$.

The same is true for asymptotic spaces: in all non-trivial cases we obtain not a single space $\mathcal{X}_0$, but rather a lot of spaces, corresponding to different nonstandard models $M$ and infinitesimals $\varepsilon$. They all correspond, in a sense, to different ways of limit passage $\varepsilon \to 0$ in the initial, na\"{i}ve definition, and all of them capture some features of the limit behavior of $^*\mathcal{X}$ at large distances.

Let us consider some examples.

\textbf{Example 1.} $\mathcal{X} = \mathbf{R}^n$ with the Euclidean distance. Then $^*\mathcal{X} = {^*\mathbf{R}^n}$, with the distance given by the same formula; so, $^*d(x,y) \in {^*\mathbf{R}}$, if $x,y \in {^*\mathbf{R}^n}$. The space $^*\mathcal{X}_\varepsilon$ is isometric to $^*\mathcal{X}$ for every $\varepsilon \in {^*\mathbf{R}_+}$. The standard part of $^*\mathbf{R}^n$ is, obviously, $\mathbf{R}^n$; thus, $\mathcal{X}_0 = \mathcal{X} = \mathbf{R}^n$. In this case the result does not depend on the model $M$ and the infinitesimal $\varepsilon$.

\textbf{Example 2.} $\mathcal{X} = (\mathbf{Z}, d)$, where $d(x,y) = |x - y|$ for every $x, y \in \mathbf{Z}$. Here $^*\mathcal{X} = (^*\mathbf{Z}, ^*d)$, where $^*\mathbf{Z}$ is the set of $^*$-integer nonstandard numbers, and the distance function $^*d$ is given by the same formula (this function assumes non-negative nonstandard integer values). Now, $^*\mathcal{X}_\varepsilon = (^*\mathbf{Z}, \varepsilon \cdot {^*d})$. Let $O$ be the zero (both in $\mathbf{Z}$ and in $^*\mathbf{Z}$). Then the finite part of $^*\mathcal{X}_\varepsilon$ is $Y = \{z \in {^*\mathbf{Z}}; |\varepsilon \cdot z| \text{ is finite}\}$.

Let us prove that $\mathrm{st}(Y) = \mathbf{R}$. For every $x \in {^*\mathbf{R}}$, let $n(x) = \left[\frac{x}{\varepsilon}\right]$ be the integer part of $\frac{x}{\varepsilon}$ ($n(x)$ is a non-standard integer). If $x, y \in \mathbf{R}$, then

\begin{align*}
^*d_\varepsilon(n(x), n(y)) \; &= \; \varepsilon|[x/\varepsilon] - [y/\varepsilon]| \\
&= |x - y + \varepsilon([x/\varepsilon] - x/\varepsilon) - \varepsilon([y/\varepsilon] - y/\varepsilon)| \tag*{(1.1)} \\
&= |x - y| + r, \qquad |r| < \varepsilon.
\end{align*}

But $\varepsilon$ is infinitesimal, so $\mathrm{st}\, ^*d_\varepsilon(n(x),n(y)) = |x - y|$. The equivalence classes consist of all $n \in {^*\mathbf{Z}}$, such that $\varepsilon \cdot n$ is contained in a \textit{monad} $\mu(x)$

of finite real number $x$, i.e. in the set of all $y \in {^*\mathbf{R}}$, such that $|x - y|$ is infinitesimally small. If $n \in Y$, then $\varepsilon n$ is in a monad $\mu(x)$ of some $x \in \mathbf{R}$. Thus we have proved that if $\mathcal{X} = (\mathbf{Z}, d)$, then $\mathcal{X}_0 = (\mathbf{R}, d)$.

\textbf{Example 3.} Let $\mathcal{X} = (X, d)$, where $X = \{2^n, n \in \mathbf{Z}_+\}$, $d$ is the usual distance in $\mathbf{R}$. In this case $^*X = \{2^n, n \in {^*\mathbf{Z}_+}\}$; let $\varepsilon \in {^*\mathbf{R}_+}$ be infinitesimal. Then we may find $n \in {^*\mathbf{Z}_+}$, such that $1/2 \leq \varepsilon \cdot 2^n < 1$ (this is true in $\mathbf{R}$, thus by the transfer principle this is true in $^*\mathbf{R}$). Let $c = \mathrm{st}(\varepsilon 2^n)$; if $m \in {^*\mathbf{Z}_+}$, $|m - n| < \infty$, then $\mathrm{st}(\varepsilon \cdot 2^m) = c \cdot 2^{m-n}$. If $m < n$, and $|m - n|$ is infinite, then $\varepsilon \cdot 2^n$ is infinitesimal, and its standard part is 0. If $m > n$, and $m - n$ is infinite, then $2^m$ is not in $Y$. Thus, $\mathcal{X}_0 = (\{0\} \cup \{c \cdot 2^n\}, d)$, where $n \in \mathbf{Z}$, $1/2 \leq c < 1$, and $c$ depends on the infinitesimal $\varepsilon$.

\textbf{Example 4.} $\mathcal{X} = (X, d)$, where $X = \{2^{2^n}, n \in \mathcal{X}_+\}$; $d(p,q) = |p - q|$. In this case $^*\mathcal{X} = (^*X, ^*d)$, where $^*X = \{2^{2^n}\}, n \in {^*\mathbf{Z}_+}$; $^*d(p,q) = |p - q|$ is a $^*\mathbf{Z}_+$-valued distance function. If $\varepsilon$ is infinitesimal, then there are the following two possibilities:

(1)There exists $n \in {^*\mathbf{Z}_+}$, such that $\varepsilon \cdot 2^{2^n}$ is finite and non-infinitesimal. Let $c = \mathrm{st}(\varepsilon \cdot 2^{2^n})$. If $m < n$, then $\varepsilon \cdot 2^{2^m} = \varepsilon \cdot 2^{2^n} \cdot 2^{(2^m-2^n)} < \varepsilon \cdot 2^{2^n} \cdot 2^{-2^{n-1}}$ is infinitesimal, while if $m > n$, then $\varepsilon \cdot 2^{2^m} = \varepsilon \cdot 2^{2^n} \cdot 2^{(2^m-2^n)} > \varepsilon \cdot 2^{2^n} \cdot 2^{2^n}$ is infinite; thus, $X_0 = \{0, c\}$.

(2)For every $n \in {^*\mathbf{Z}_+}$, $\varepsilon \cdot 2^{2^n}$ is either infinitesimal (e.g. for $n$ finite), or infinite (e.g. for $n > 1/\varepsilon^2$). In this case, $X_0 = \{0\}$.

We see that in the examples 3 and 4 the asymptotic space $\mathcal{X}_0$ depends on an infinitesimal $\varepsilon$. In the forthcoming sections we shall see that asymptotic space may depend dramatically on the nonstandard model, too.

\textbf{4.} The next example, which is the subject of this work, is the asymptotic space of the Lobachevsky (or hyperbolic) plane $\mathcal{H}$. This is the simplest representative of the class of hyperbolic spaces. The general principle is, that for every hyperbolic space $\mathcal{X}$ its asymptotic space $\mathcal{X}_0$ is an $\mathbf{R}$-tree; this means that for every two points $x, y \in X_0$ there exists unique geodesic segment $[x, y]$, joining $x$ and $y$, and if $[x, y]$ and $[y, z]$ are two geodesic segments in $X_0$, having only one common point $y$, then their union is a geodesic segment $[x, z]$, connecting $x$ and $z$ (i.e. there is no cycles)(see [5]).

The problem is, how to describe such tree-like space for concrete hyperbolic space $\mathcal{X}$. If $\mathcal{X} = \mathcal{H}$, then we know that it is a homogeneous space. It is natural to conjecture that its asymptotic space $\mathcal{H}_0$ is also homogeneous. Thus, $\mathcal{H}_0$ is a tree, branching at every point. But this is not a satisfactory definition; to describe a metric space, we should construct the set of points, and the distance function. The aim of this work is to give an explicit description of the space $\mathcal{H}_0$. We are going to construct a (standard) metric space, which is isometric to $\mathcal{H}_0$. The success of this business is not complete and depends on the nonstandard model $M$ of the superstructure $\mathcal{S}$. For every model $M$ we define an asymptotic space $\mathcal{H}_{0,M}$, and construct a metric space $\mathcal{F}_M$, such that $\mathcal{H}_{0,M}$ may be isometrically embedded in the space $\mathcal{F}_M$. On the other hand, for every model $M$ we construct some other model $M'$, such that $\mathcal{F}_M$ may be isometrically embedded in $\mathcal{H}_{0,M'}$. And for some distinguished models, called \textit{saturated} models, we prove that the spaces $\mathcal{F}_M$ and $\mathcal{H}_{0,M}$ \textit{coincide}.

To give the idea of the definition of spaces $\mathcal{F}_M$, we describe simpler spaces $\mathcal{C}$ and $\mathcal{D}$, which are isometric to some subcones of $\mathcal{H}$. (see [8]). The space $\mathcal{C}$ is the set of pairs $(f, a)$, where $a \in \mathbf{R}, a \geq 0$, and $f = f(x)$ is a continuous function, defined on the segment $0 \leq x \leq a$, and such that $f(0) = 0$. If $\alpha = (f, a)$, $\beta = (g, b)$ are two elements of $\mathcal{C}$, then we define the distance

\begin{equation*}
d_{\mathcal{C}}(\alpha, \beta) = (b - c) + (a - c), \tag{1.2}
\end{equation*}

\noindent where

\begin{equation*}
c = \sup\{s \mid f(x) = g(x), \quad 0 \leq x \leq s\} \tag{1.3}
\end{equation*}

\noindent ($c$ is called the moment of separation of $f$ and $g$). It is easy to see that $d_{\mathcal{C}}$ defines a true metric in $\mathcal{C}$, and, as it was proved in [7], $C$ with this metric is a homogeneous $\mathbf{R}$-tree.

The space $\mathcal{D}$ is also a set of pairs $(f, a)$, where $a \geq 0$, $f(x)$ is defined for $0 \leq x \leq a$, and $f(x) = 0$ everywhere but at a finite set of points. The distance is defined by the same formulas as for the space $\mathcal{C}$. We have proved in [7] that both $\mathcal{C}$ and $\mathcal{D}$ are asymptotic subcones of $\mathcal{H}$, and these spaces are non-isometric.

We may define more general spaces as follows. Let $\Lambda$ be a linearly ordered set with order relation $<$. Let $g(\lambda)$ be a real-valued function on $\Lambda$, which is non-decreasing with respect to the order on $\Lambda$, and whose values cover completely some segment $(g_1, g_2)$ (finite, infinite, or semi-infinite) in $\mathbf{R}$.

Our space $\mathcal{F}$ consists of the pairs $(f(\lambda), \bar{\lambda})$, where $\bar{\lambda} \in \Lambda$, and $f(\lambda)$ is a real function (of some specified class depending on the model $M$ and infinitesimal $\varepsilon$) defined for $\lambda_0 \leq \lambda \leq \bar{\lambda}$, where $\lambda_0$ is some fixed element of  $\Lambda$. If $\alpha_1 = (f_1, \bar{\lambda}_1)$, $\alpha_2 = (f_2, \bar{\lambda}_2)$ are two elements of the space $\mathcal{F}$, then we define

\begin{equation*}
d_{\mathcal{F}}(\alpha_1, \alpha_2) = (g(\bar{\lambda}_1) - c) + (g(\bar{\lambda}_2) - c), \tag{1.4}
\end{equation*}

\noindent where

\begin{equation*}
c = \sup\{s \mid f_1(\lambda) = f_2(\lambda) \text{ for all } \lambda, \text{ such that } g(\lambda) \leq s\}. \tag{1.5}
\end{equation*}

\medskip
\noindent For the same initial space $\mathcal{H}$ there is a plenty of spaces $\mathcal{F}_M$, corresponding to different models $M$ (in particular, there exist such spaces having arbitrarily big cardinality).

\bigskip
\noindent \textbf{5.} It should be noted that our constructions resemble the works of G. W. Brumfiel [1] devoted to the study of hyperbolic planes over nonarchimedean valued field $F$. Namely, for some classes of these fields (called ``microbial'') it is possible to define a pseudometric on the hyperbolic plane $HF^2$ over $F$, using the values in $F$ or in its algebraic extension. The main result of G. W. Brumfiel is that the metric space quotient $TF^2$ of $HF^2$ is an $R$-tree. This result is further used in the study of compactifications of the Teichm\"{u}ller space [2]. Note that ${}^{*}\mathbf{R}$ is not a valued field (it is too big), and our problem and results are different from that of G. W. Brumfiel. (I am thankful to Pierre Pansu for pointing me on the works of G. W. Brumfiel.)

\bigskip
\noindent \textbf{6.} Our work has the following structure. In Section 2 we collect necessary formulas from the analytic geometry of the Lobachevsky plane.

\medskip
\noindent In Section 3 we develop the necessary tools of the nonstandard number theory. The nonstandard numbers form a module ${}^{*}\mathbf{R}$ over the field $\mathbf{R}$. It has infinite dimension, and we define a ``basis'', such that every nonstandard number has unique ``decomposition'' relative to this basis with standard real coefficients. It turns out that the problem of ``synthesis'' of non-standard number cannot, in general, be solved in the given nonstandard $M$, and requires a larger model $M'$. The models where the synthesis is always possible are just the \textit{saturated} (i.e. sufficiently large) ones [3]. Existence of these models depends on the Generalized Continuum Hypothesis (GCH).

\medskip
\noindent In Section 4 we develop the nonstandard geometry of the Lobachevsky plane; we use explicit formulas of the analytic non-Euclidean geometry, and the simple computations, similar to those made in the above examples, give us explicit description of the asymptotic space. We construct the space $\mathcal{F}_M$ such that the asymptotic space $\mathcal{H}_{0,M}$ is isometrically embedded in $\mathcal{F}_M$. On the other hand, we construct here another, larger model $M'$, such that $\mathcal{F}_M$ may be isometrically embedded into $\mathcal{H}_{0,M'}$.

We prove in Section 5 that if the model $M$ is saturated, than $\mathcal{H}_{0,M}$ coincides with $\mathcal{F}_M$, and thus in this case we have an adequate description of $\mathcal{H}_{0,M}$.

At last, we establish in Section 6 the $\mathbf{R}$-tree structure of the space $\mathcal{F}_M$. While there is a general principle that asymptotic space of a hyperbolic space is an $\mathbf{R}$-tree, it is interesting to see how this structure materializes for concrete space. The asymptotic space $\mathcal{H}_{0,M}$ is always homogeneous; but it is isometric to the space $\mathcal{F}_M$ only if $M$ is a saturated model. Otherwise homogeneity of the space $\mathcal{F}$ is unclear.

\noindent\textbf{7.} For every unbounded metric space $\mathcal{X}$ we may ask, what is its asymptotic space. The next spaces of great interest are the groups of volume preserving diffeomorphisms with $L^2$-metric, and the groups of symplectomorphisms with the Hofer metric. These groups are important in the fluid dynamics and in the Hamiltonian dynamics, respectively, and in particular their asymptotic structure is very interesting. We may formulate here one conjecture. Let us call the metric space \textit{large}, if its asymptotic space in some (or every?) nonstandard model $M$ contains a subspace isometric to $\mathcal{F}_M$.

\medskip
\noindent\textbf{Conjecture.} The group of area preserving diffeomorphisms of 2-dimensional domain with $L^2$-metric, and the group of symplectomorphisms with the Hofer metric are large spaces.

\medskip
\noindent\textbf{Acknowledgements.} I am grateful to Misha Shubin for the first acquaintance with the nonstandard analysis, and to Micha Gromov for fruitful discussions. I am especially thankful to my colleagues Yoram Hirshfeld and Moti Gitik for their interest to this work and a lot of useful discussions and improvements.

Most of this work had been done in the Institut des Hautes \'{E}tudes Scientifiques. I am thankful to this institution for its hospitality and excellent work conditions.

\section{Analytic geometry in the Lobachevsky plane}

\noindent The Lobachevsky plane $\mathcal{H} = (H,d)$ may be realized as a unit disk $|z| < 1$ in the complex $z$-plane with the distance function

\begin{equation*}
d(z_1, z_2) = \frac{1}{2} \log \frac{1 + A}{1 - A}, \tag{2.1}
\end{equation*}

\noindent where

\begin{equation*}
A = \left| \frac{z_1 - z_2}{z_1 \bar{z}_2 - 1} \right|. \tag{2.1}
\end{equation*}

Let us introduce polar coordinates in $H$: $x = (\rho,\varphi)$, $0 \leq \rho < \infty$, $0 \leq \varphi < 2\pi$, where $\rho$ is the Lobachevsky distance between 0 and $x$, and $\varphi$ is the polar angle. If $x_1 = (\rho_1,\varphi_1)$, $x_2 = (\rho_2,\varphi_2)$, then

\begin{equation*}
d(x_1, x_2) = \frac{1}{2} \log \frac{1 + A}{1 - A}, \tag{2.3}
\end{equation*}

\noindent where

\begin{equation*}
A^2 = 1 - \frac{8}{(2 - \varepsilon^2)(t + 1/t)^2 + \varepsilon^2(s + 1/s)^2}; \tag{2.4}
\end{equation*}

\begin{equation*}
s^2 = e^{\rho_1 + \rho_2}; \quad t^2 = e^{\rho_1 - \rho_2}; \quad \varepsilon^2 = 1 - \cos(\varphi_1 - \varphi_2). \tag{2.5}
\end{equation*}

Simple asymptotic analysis of these formulas shows, that if $\rho_1 = N \cdot R_1$, $\rho_2 = N \cdot R_2$, $\varphi_1 - \varphi_2 = e^{-N\Phi}$, then for $N \to \infty$,

\begin{equation*}
A^2 = 1 - \frac{8}{B}, \tag{2.6}
\end{equation*}

\noindent where

\begin{align*}
B \; = \; & (2 - e^{-2N\Phi})(e^{N(R_1-R_2)} + 2 + e^{N(R_1-R_2)}) \\
& +e^{-2N\Phi}(e^{N(R_1+R_2)} + 2 + e^{-N(R_1+R_2)}); \tag{2.7}
\end{align*}

\noindent If $R_1 > R_2 > 0$, $\Phi > 0$, then

\begin{equation*}
A^2 = 1 - \frac{8}{2e^{N(R_1-R_2)} + e^{-2N\Phi} \cdot e^{N(R_1+R_2)} + O(e^{-NQ})}, \tag{2.8}
\end{equation*}

\thispagestyle{plain}

\begin{equation*}
Q = \min(0,\; 2\Phi - R_1 + R_2). \tag{2.9}
\end{equation*}

\noindent Thus,

\begin{align*}
\frac{1}{2} \log \frac{1+A}{1-A} \; = \; & \log \left[ e^{N|R_1-R_2|} + e^{N(R_1+R_2-2\Phi)} \right. \\
& \left. + O(1) + O\!\left(e^{N(|R_1-R_2|-2\Phi)}\right) \right] + O(1); \tag{2.10}
\end{align*}

\begin{align*}
\frac{1}{N} d(x_1,x_2) \; = \; & \frac{1}{N} \log \left[ e^{N|R_1-R_2|} + e^{N(R_1+R_2-2\Phi)} + O(1) \right. \\
& \left. + O\!\left(e^{N(|R_1-R_2|-2\Phi)}\right) \right] + O\!\left(\frac{1}{N}\right). \tag{2.11}
\end{align*}

\noindent So,

\begin{equation*}
\frac{1}{N} d(x_1,x_2) = \max\!\left(|R_1-R_2|,\, R_1+R_2-2\Phi\right) + O\!\left(\frac{1}{N}\right). \tag{2.12}
\end{equation*}

\section{Decomposition of a nonstandard number}

\noindent In order to describe the asymptotic space of the hyperbolic plane, we have to introduce some adequate techniques. Consider a nonstandard extension ${}^{*}\mathbf{R}$ of the field $\mathbf{R}$ of real numbers. Since $\mathbf{R} \subset {}^{*}\mathbf{R}$, ${}^{*}\mathbf{R}$ is an $\mathbf{R}$-module. Our goal here is to find a ``basis'' in ${}^{*}\mathbf{R}$, so that every nonstandard number could be represented as a linear combination of elements of this basis with standard real coefficients, and this representation be unique.

Our construction requires a lot of use of the Axiom of Choice (AC). We may unify all the choices, using the following trick (proposed by Yoram Hirshfeld). Let $M$ be a model of the superstructure $\mathcal{S}$ (based on the set $\mathbf{R}$), and ${}^{*}\mathbf{R}$ the set of nonstandard numbers in the model $M$. Let us fix some one-to-one correspondence between ${}^{*}\mathbf{R}$ and its cardinal $|{{}^{*}\mathbf{R}}|$; let $\varphi(x)$ denote the ordinal corresponding to $x \in {}^{*}\mathbf{R}$. By the definition of cardinal as the least ordinal having power $|{{}^{*}\mathbf{R}}|$, $\varphi(x) < |{{}^{*}\mathbf{R}}|$ for all $x \in {}^{*}\mathbf{R}$. This means that we have introduced some well-ordering in the set ${}^{*}\mathbf{R}$. Let us denote this relation by $\prec$: $x \prec y$, if $\varphi(x) < \varphi(y)$. From now on, if we have to pick an element from some set $B \subset {}^*\mathbf{R}$, we shall take the minimal element of $B$ with respect to the order $\prec$.

Let $x, y \in {}^*\mathbf{R}$, $x, y \neq 0$; we write $x \sim y$, if both $x/y$ and $y/x$ are finite elements of ${}^*\mathbf{R}$ (i.e. if $|\frac{x}{y}| < K$, $|\frac{y}{x}| < K$ for some $K \in \mathbf{R}$). $\sim$ is an equivalence relation; the factor-set of ${}^*\mathbf{R}$ with respect to this relation is denoted by $\Lambda$. This is an ordered abelian group. For every $x \in {}^*\mathbf{R}$, we denote by $\lambda(x)$ the image of $x$ in $\Lambda$. Then the group operation in $\Lambda$ is defined as $\lambda(x) + \lambda(y) = \lambda(x \cdot y)$. The order in $\Lambda$ is induced by the natural order in ${}^*\mathbf{R}$: $\lambda(x) < \lambda(y)$, if $\lambda(x) \neq \lambda(y)$, and $0 < |x| < |y|$. We denote this relation by the same symbol $<$.

For every $\lambda \in \Lambda$, let us define its representative $a_\lambda \in \lambda$ as the minimal element of $\lambda$ with respect to the order relation $\prec$ in ${}^*\mathbf{R}$.

Let $x \in {}^*\mathbf{R}$; we shall define sequences $x_i \in {}^*\mathbf{R}, \lambda_i \in \Lambda, c_i \in \mathbf{R}, s_i \in {}^*\mathbf{R}$ in the following way:

\begin{align*}
&x_0 = x; \quad \lambda_0 = \lambda(x_0); \quad c_0 = \operatorname{st}\!\left(\frac{x_0}{a_{\lambda_0}}\right); \quad s_0 = 0; \\
&x_1 = x_0 - c_0 a_{\lambda_0}; \quad \lambda_1 = \lambda(x_1); \quad c_1 = \operatorname{st}\!\left(\frac{x_1}{a_{\lambda_1}}\right); \quad s_1 = s_0 + |{c_0 \cdot a_{\lambda_0}}|; \\
&\cdots\cdots\cdots\cdots\cdots\cdots\cdots\cdots\cdots\cdots\cdots\cdots\cdots\cdots\cdots\cdots\cdots\cdots \\
&x_{i+1} = x_i - c_i a_{\lambda_i}; \quad \lambda_{i+1} = \lambda(x_{i+1}); \quad c_{i+1} = \operatorname{st}\!\left(\frac{x_{i+1}}{a_{\lambda_{i+1}}}\right); \\
&s_{i+1} = s_i + |{c_{i} \cdot a_{\lambda_{i}}}|; \tag{3.1}\\
&\cdots\cdots\cdots\cdots\cdots\cdots\cdots\cdots\cdots\cdots\cdots\cdots\cdots\cdots\cdots\cdots\cdots\cdots\cdots\cdots
\end{align*}

Here $i$ is natural, and $\operatorname{st}(a)$ means the standard part of a finite nonstandard number $a$. We see that either at some step we obtain $x_{i+1} = 0$, and the decomposition of $x$ into a linear combination of $a_{\lambda_i}$ is over, or we may go on. At every non-terminal step we obtain $c_i \neq 0$, and $\lambda_i > \lambda_{i+1}$ in the sense of the order in $\Lambda$ induced by the usual order in ${}^*\mathbf{R}$.

After we have constructed $x_i, c_i, \lambda_i, s_i$ for all $i < \omega$, we may define $x_\omega, c_\omega, \lambda_\omega, s_\omega$ in the following way. Consider the set $X_{(c_i, \lambda_i), i<\omega}$ of all $y \in {}^*\mathbf{R}$, such that the above process gives the sequence $(c_i)$ of real coefficients and $(\lambda_i)$ of classes for all $i < \omega$. This set is not empty; it contains $x$ itself. Let us choose the \textit{minimal} element of $X_{(c_i, \lambda_i), i<\omega}$ with respect to the order $\prec$; denote it by $a_{(c_i, \lambda_i), i<\omega}$. Now we may define

\begin{equation*}
s_\omega = a_{(c_i, \lambda_i),\, i<\omega};\quad x_\omega = x - s_\omega;\quad \lambda_\omega = \lambda(x_\omega);\quad c_\omega = \operatorname{st}\!\left(\frac{x_\omega}{a_{\lambda_\omega}}\right), \tag{3.2}
\end{equation*}

\noindent and then define $x_{\omega+1} = x_\omega - c_\omega a_{\lambda(x_\omega)}$, etc. for all ordinals $i < 2\omega$.

In general, we define the sequences $(x_i),(c_i),(\lambda_i),(s_i)$ for all ordinals $i$ using the transfinite recursion. Suppose that $\alpha$ is an ordinal, and for all ordinals $\beta < \alpha$ we have already defined $x_\beta, c_\beta$. If $\alpha = \beta + 1$ for some $\beta$, then we define simply $x_\alpha = x_\beta - c_\beta a_{\lambda_\beta}$, etc. If $\alpha$ is a limit ordinal, then denote by $X_{(c_\beta, \lambda_\beta), \beta < \alpha}$ the set of all nonstandard numbers such that the above process delivers the sequence $(c_\beta, \lambda_\beta)$ for all $\beta < \alpha$. This set is not empty for it contains $x$ itself (and all nonstandard numbers $z$ such that $z - x = o(a_{\lambda_\beta})$ for all $\beta < \alpha$). Let $a_{(c_\beta, \lambda_\beta), \beta < \alpha}$ be the minimal element of the set $X_{(c_\beta, \lambda_\beta), \beta < \alpha}$ with respect to the order relation $\prec$. Then define

\begin{align*}
s_\alpha &= a_{(c_\beta, \lambda_\beta),\, \beta < \alpha};\\
x_\alpha &= x - a_{(c_\beta, \lambda_\beta),\, \beta < \alpha};\\
\lambda_\alpha &= \lambda(x_\alpha); \quad c_\alpha = \operatorname{st}\!\left(\frac{x_\alpha}{a_{\lambda_\alpha}}\right). \tag{3.3}
\end{align*}

Thus we have defined $x_\alpha$ for all ordinals $\alpha$. It is possible that $x_\alpha = 0$ for some $\alpha$; in this case we interrupt the process.

Note that at every step we go strictly down the ordered set $\Lambda$: $\lambda_\alpha < \lambda_\beta$, if $\alpha > \beta$. So, the number of steps, i.e. $|\alpha|$, cannot exceed $|\Lambda|$. This means that for every $x \in {}^*\mathbf{R}$ the process terminates at some step, i.e. $x_\alpha = 0$ for some $\alpha, |\alpha| \leq |\Lambda|$.

The sequence $(\lambda_i)$, associated with the number $x \in {}^*\mathbf{R}$, is the \textit{spectrum} of the nonstandard number $x$, and is denoted by $\operatorname{spec}(x)$. The spectrum depends implicitly on the choice of representatives $a_\lambda$, or $a_{(c_i, \lambda_i)}$, $i < \alpha$, playing the role of a ``basis'' in ${}^*\mathbf{R}$ over $\mathbf{R}$. If $\operatorname{spec}(x)$ has the last term with index $\alpha$, then $x = s_\alpha$, and if $x_i$ is defined for all $i < \alpha$, $\alpha$ a limit ordinal, then $x = a_{(c_i, \lambda_i)}$, $i < \alpha$.

\textbf{Lemma 3.1} \textit{Nonstandard number $x \in {}^*\mathbf{R}$ is uniquely characterized by the sequence $(c_i, \lambda_i)$.}

\textbf{Proof.} Suppose that $x, y \in {}^*\mathbf{R}$, and $x \neq y$. Let $(c_i, \lambda_i), (d_i, \mu_i)$ be the sequences, corresponding to $x$ and $y$. By the transfinite induction we prove that if $\lambda_i = \mu_i$, $c_i = d_i$ for all $i \leq \alpha$, then $|x - y| < a_{\lambda_\beta}$. If $\alpha$ is the first ordinal such that $c_\alpha = d_\alpha = 0$, then $x = y = s_\alpha$; if $\alpha$ is the first limit ordinal, such that $c_\alpha = d_\alpha = 0$, then $x = y = a_{(c_i, \lambda_i)}$.

\begin{flushright}
Q.E.D.
\end{flushright}

\bigskip

\noindent Important property of the spectra is stated by the following

\bigskip
\textbf{Theorem 3.1} $|\operatorname{spec}(x)| < |{}^*\mathbf{R}|$ \textit{for every nonstandard number $x$.}

\bigskip
\textbf{Proof.} Previously we have fixed a one-to-one correspondence between the elements of the set ${}^*\mathbf{R}$, on one hand, and all ordinals less than $|{}^*\mathbf{R}|$, on the other. Let $\varphi(x)$ denote the ordinal, corresponding to $x \in {}^*\mathbf{R}$.

If $\beta, \gamma$ are two limit ordinals, and $\beta < \gamma < \alpha$, where the ordinal $\alpha$ is the length of $\operatorname{spec}(x)$, then $X_{(c_i, \lambda_i), i<\beta} \supset X_{(c_i, \lambda_i), i<\gamma}$. Therefore $a_{(c_i, \lambda_i), i<\beta} \prec a_{(c_i, \lambda_i), i<\gamma}$, because we choose always the \textit{minimal} element in every set $X_{(c_i, \lambda_i), i<i_0}$ with respect to the order relation $\prec$, for every limit ordinal $i_0$. On the other hand, $a_{(c_i, \lambda_i), i<i_0} \prec x$ for every $i_0$, and in the case of equality the decomposition process is over. So, the cardinality of the set of limit ordinals, appearing as indices in $\operatorname{spec}(x)$, is less than $\varphi(x)$. On the other hand, $|\varphi(x)| < |{}^*\mathbf{R}|$, because $\varphi({}^*\mathbf{R})$, by our construction, is a \textit{cardinal}, i.e. the minimal ordinal with cardinality $|{}^*\mathbf{R}|$, and the cardinality of limit ordinals, not exceeding $\varphi(x)$, is equal to $|\varphi(x)|$ if $\varphi(x) \geq \omega \cdot \omega$; otherwise $|\varphi(x)| \leq \aleph_0$. Thus, $|\operatorname{spec}(x)| < |{}^*\mathbf{R}|$.

\begin{flushright}
Q.E.D.
\end{flushright}

\medskip
(I am thankful to Prof. Yoram Hirshfeld for the idea of this theorem.)

Given $x \in {}^*\mathbf{R}$, let us define a function

\begin{equation*}
f_x(\lambda) = \left\{ \begin{array}{ll} c_i, & \text{if } \lambda = \lambda_i \in \operatorname{spec} x; \\ 0, & \text{if } \lambda \notin \operatorname{spec} x. \end{array} \right. \tag{3.4}
\end{equation*}

This is a real-valued function; it is different from 0 on the set $\operatorname{spec}x$, which is well ordered with respect to the order in $\Lambda$, inverse to the natural order; the cardinality of the set $\operatorname{spec}(x)$ is always less than $|{}^*\mathbf{R}|$. Lemma 3.1 shows that if $x, y \in {}^*\mathbf{R}, x \neq y$, then $f_x(\lambda) \neq f_y(\lambda)$. But it is, in general, unclear, that
for \textit{arbitrary} function $f(\lambda)$, satisfying the previous conditions, there exists an $x \in {}^*\mathbf{R}$, such that $f(\lambda) \equiv f_x(\lambda)$. This question is important for us, and we shall return to it below.

Now let us consider the similar construction for the nonstandard circle ${}^*S^1$. The nonstandard semi-interval ${}^*[0,2\pi)$ contains points which are infinitesimally close to $2\pi$; so, we shall regard ${}^*S^1$ as ${}^*[0,2\pi]$ with identified points 0 and $2\pi$.

The first step of our construction yields a number $a_{\lambda_0}$, where $\lambda_0$ is the set of finite and non-infinitesimal numbers. Let us set $a_{\lambda_0} = 1$; this changes nothing in our analysis.

If $x \in {}^*S^1$, then $c_0 = \mathrm{st}\!\left(\frac{x}{a_{\lambda_0}}\right) = \frac{\mathrm{st}(x)}{\mathrm{st}(a_{\lambda_0})}$. It is evident, that $0 \leq c_0 \leq 2\pi$. If $\mathrm{st}\, x = 0$ or $\mathrm{st}\, x = 2\pi$, we define $c_0 = 0$. Then $c_0$ may be regarded as an element of standard $S^1$. For the next step, if $c_0 \neq 0$, we proceed as before.

So, if $x \in {}^*S^1$, then the function $f_x(\lambda)$ is defined on $\Lambda$; it is zero for $\lambda > \lambda_0$, $f_x(\lambda_0) \in S^1$, and $f_x(\lambda) \in \mathbf{R}$ for $\lambda < \lambda_0$, with all properties as above.

\section*{4 \quad Asymptotic space}

Let $M$ be a nonstandard model of the universum $\mathcal{S}$ (i.e. $M$ contains ${}^*\mathbf{R}, {}^*\mathcal{P}(\mathbf{R})$, etc.; in particular, every real-valued function may be extended to the nonstandard one). Let ${}^*\mathcal{H} = ({}^*H, {}^*d)$ be a nonstandard hyperbolic plane, where ${}^*H$ is the set of pairs $(\rho,\varphi)$, $\rho \in {}^*\mathbf{R}$, $\varphi \in {}^*S^1$, and the distance ${}^*d$ is defined by the same formulas (2.3)-(2.5), where the functions $\exp$, $\log$ are understood as the nonstandard extensions of the standard ones.

Let $\varepsilon \in {}^*\mathbf{R}$ be an infinitesimal (super)real number. Consider the space ${}^*\mathcal{H}_{\varepsilon} = ({}^*H, {}^*d_{\varepsilon})$, where $d_{\varepsilon}(x,y) = \varepsilon \cdot d(x,y)$, $x,y \in {}^*H$. Let $Y \subset {}^*H, Y = \{(\rho,\varphi) \in {}^*H \mid \varepsilon \cdot \rho \text{ is finite}\}$. The distances between the points of $Y$ are finite (and nonstandard). Let $x,y \in Y$; we say that $x \approx y$, if $d_{\varepsilon}(x,y)$ is infinitesimal. The space of the equivalence classes, with the distance equal to the standard part of $d(x,y)$ between the representatives, is a standard metric space, which by definition is called the \textit{asymptotic space} of $\mathcal{H}$ and denoted by $\mathcal{H}_{0,M}$ (we write the subscript $M$, because this space depends on the model). Our goal is to describe its structure.

If $\varepsilon \in {}^*\mathbf{R}$ is positive and infinitesimal, $x_1, x_2 \in {}^*\mathcal{H}_{\varepsilon}$, $x_1 = (\rho_1, \varphi_1)$, $x_2 = (\rho_2, \varphi_2)$, then

\begin{equation*}
\operatorname{st} {}^*d_{\varepsilon}(x_1, x_2) = \operatorname{st}(\varepsilon \log(e^{|\rho_1 - \rho_2|} + \beta^2 e^{\rho_1 + \rho_2})), \tag{4.1}
\end{equation*}

\begin{equation*}
\beta^2 = 1 - \cos(\varphi_1 - \varphi_2). \tag{4.2}
\end{equation*}

If $\varphi_1 - \varphi_2 = e^{-\Phi}$, then

\begin{equation*}
\operatorname{st} {}^*d_{\varepsilon}(x_1, x_2) = \operatorname{st} \max(\varepsilon|\rho_1 - \rho_2|,\, \varepsilon(\rho_1 + \rho_2 - 2\Phi)). \tag{4.3}
\end{equation*}

This expression may be associated with the function $f_{\varphi}(\lambda)$, defined by the formula (3.4), where $\varphi$ is substituted instead of $x$.

We define another function on $\Lambda$,

\begin{equation*}
g(\lambda) = \operatorname{st}(-\varepsilon \log a_{\lambda}). \tag{4.4}
\end{equation*}

For any two real-valued functions on $\Lambda$, $f_1(\lambda), f_2(\lambda)$, defined for $\lambda_i \leq \lambda \leq \lambda_0$ ($i = 1,2$), let us define the distance

\begin{equation*}
\delta((f_1, \lambda_1),(f_2, \lambda_2)) = (g(\lambda_1) - \tilde{g}) + (g(\lambda_2) - \tilde{g}), \tag{4.5}
\end{equation*}

where

\begin{equation*}
\tilde{g} = \sup\{h \mid f_1(\lambda) = f_2(\lambda) \text{ for all } \lambda \leq \lambda_0, \text{ such that } g(\lambda) < h\}. \tag{4.6}
\end{equation*}

\textbf{Theorem 4.1} \textit{Let $x_i = (\rho_i, \varphi_i) \in Y, i = 1,2$. Set $\lambda_i = \lambda(e^{-\rho_i})$. Then}

\begin{equation*}
\operatorname{st}\, d_{\varepsilon}(x_1, x_2) = \delta((f_{\varphi_1}, \lambda_1),(f_{\varphi_2}, \lambda_2)). \tag{4.7}
\end{equation*}

\textbf{Proof.} Consider the functions $f_{\varphi_1}(\lambda), f_{\varphi_2}(\lambda)$. They are supported on the sets $\operatorname{spec}(\varphi_1), \operatorname{spec}(\varphi_2)$; these sets and their union are well ordered with respect to the inverse order on $\Lambda$.

Consider the set $\Sigma \subset \operatorname{spec}(\varphi_1) \cup \operatorname{spec}(\varphi_2)$,

\begin{equation*}
\Sigma = \left\{ \lambda \in \Lambda \,\Big|\, |f_{\varphi_1}(\lambda)| + |f_{\varphi_2}(\lambda)| > 0, \text{ and } f_{\varphi_1}(\lambda) \neq f_{\varphi_2}(\lambda) \right\}. \tag{4.8}
\end{equation*}

There exists a maximal in $\Lambda$ element $\tilde{\lambda} \in \dot{\Sigma}$.
Thus, $\varphi_1 - \varphi_2 \in \tilde{\lambda}$, and

\begin{equation*}
\Phi = -\log|\varphi_1 - \varphi_2| = -\log a_{\tilde{\lambda}} + O(1); \tag{4.9}
\end{equation*}

thus,

\begin{equation*}
\tilde{g} = \min\!\left(\operatorname{st}\!\left(-\varepsilon \log a_{\tilde{\lambda}}\right), \lambda_1\right) = \min\!\left(\varepsilon\Phi, \varepsilon\rho_1\right), \tag{4.10}
\end{equation*}

because $\varepsilon$ is infinitesemal.

Assume that $\rho_1 \leq \rho_2$; then

\begin{equation*}
\operatorname{st} d_{\varepsilon}(x_1, x_2) = \operatorname{st} \max\!\left(\varepsilon|\rho_1 - \rho_2|, \varepsilon(\rho_1 + \rho_2 - 2\Phi).\right. \tag{4.11}
\end{equation*}

On the other hand,

\begin{equation*}
g(\lambda_i) = \operatorname{st}(-\varepsilon \log a_{\lambda_i}) = \operatorname{st}\!\left(-\varepsilon\left(\log(e^{-\rho_i} + O(1)\right) = \varepsilon\rho_i, \quad i = 1,2; \right. \tag{4.12}
\end{equation*}

\begin{equation*}
\tilde{g} = \min(\varepsilon\Phi, \varepsilon\rho_1). \tag{4.13}
\end{equation*}

Thus,

\begin{equation*}
\operatorname{st}\, d_{\varepsilon}(x_1, x_2) = \operatorname{st} \max\!\left(\varepsilon(\rho_2 - \rho_1), \varepsilon(\rho_1 + \rho_2 - 2\Phi)\right); \tag{4.14}
\end{equation*}

\begin{align*}
\delta((f_{\varphi_1}, \lambda_1),(f_{\varphi_2}, \lambda_2)) &= (\varepsilon\rho_1 - \tilde{g}) + (\varepsilon\rho_2 - \tilde{g}) \\
&= \operatorname{st} \max(\varepsilon(\rho_1 - \Phi), 0) + \operatorname{st} \max(\varepsilon(\rho_2 - \Phi), \varepsilon(\rho_2 - \rho_1)); \tag{4.15}
\end{align*}

It is easy to see that in all three possible cases $(0 \leq \Phi \leq \rho_1; \rho_1 < \Phi \leq \rho_2; \rho_2 < \Phi)$ we have st $d_{\varepsilon}(x_1,x_2) = \delta((f_{\varphi_1},\lambda_1),(f_{\varphi_2},\lambda_2))$.

\hfill QED

For a given nonstandard model $M$ of the superstructure $\mathcal{S}$ let us denote by $\mathcal{F}_M$ the space of pairs $(f, \lambda_1)$, where $\lambda_1 \leq \lambda_0$, $f(\lambda)$ is a real-valued function, defined for $\lambda_1 \leq \lambda \leq \lambda_0$, $f(\lambda_0) \in S^1$, $f(\lambda) = 0$ for all $\lambda$ except some well-ordered set $\Sigma_f$, such that $|\Sigma_f| < |{}^*\mathbf{R}|$, and the distance $\delta$ between elements of $\mathcal{F}_M$ is defined by (4.15). Then our Theorem states that the asymptotic space $\mathcal{H}_{0,M}$ is isometrically embedded into $\mathcal{F}_M$. Thus, we have found some, in general, larger space, containing the asymptotic space.

\setlength{\parindent}{20pt}
\setlength{\parskip}{0pt}

\noindent\hspace{20pt}We may try to solve, in a sense, an opposite problem and construct, for a given nonstandard model $M$ of the real universum, another, larger model $M'$, such that the space $\mathcal{F}_M$, corresponding to the model $M$, may be isometrically embedded into the space $\mathcal{H}_{0,M'}$ in the model $M'$. This means the following: for every model $M$, let us denote by $\Lambda_M$ the set $\Lambda$ of equivalence classes of nonstandard real numbers in the model $M$. If $M' = (A', E', I')$ is an elementary extension of the model $M = (A, E, I)$, i.e. $A \subset A'$, and every formula in the language $\mathcal{L}$, which is true in $M$, remains true in $M'$, then $\Lambda_M \subset \Lambda_{M'}$. So, for every function $f(\lambda)$ belonging to the above class we have to construct a number $\varphi \in {}^*\mathbf{R}_{M'}$, such that if we find the spectrum of $\varphi$, and construct the function $f_\varphi(\lambda)$, then $f_\varphi(\lambda) \equiv f(\lambda)$ on $\Lambda$, and $f_\varphi \equiv 0$ on $\Lambda_{M'} \setminus \Lambda_M$. So, let ${}^*\mathbf{R}$ be the field of nonstandard real numbers in a given model $M$, and let $\Lambda$ be the ordered set of their equivalence classes with respect to the order relation $\sim$. Suppose that $(\lambda_i)_{i<\alpha} \subset \Lambda$ is a decreasing sequence, and $(c_i)_{i<\alpha}$ is a sequence of real numbers; assume that $|\alpha| < |{}^*\mathbf{R}_M|$, and $\lambda_0 \geq \lambda \geq \bar{\lambda}_1$. We have to give some meaning to a formal expression $\sum_{i<\alpha} c_i a_{\lambda_i}$, i.e. to assign to this formal sum a number in an extended set ${}^*\mathbf{R}_{M'}$. This is similar to the problem, how to give a meaning to a formal asymptotic series $\sum_{i<\alpha} c_i f_i(x)$, where $f_i = o(f_j)$, if $i > j$, and $\alpha$ is a given ordinal. (See [7] for the analysis of this problem from the nonstandard viewpoint.)

So, suppose that $\lambda_{i}$ is a decreasing sequence in $\Lambda$, and for every ordinal $i$ (less than some $\alpha$) there is defined a real number $c_{i}$; moreover, in every class $\lambda_{i}$ there is defined a representative $a_{\lambda_i}$. We define a sequence $s_i$ by the transfinite recursion. At the same time we construct a growing sequence of models $M_{i}$.

For all finite $i$ define $s_i = s_0a_{\lambda_0} + \dots + c_ia_{\lambda_i}$, $M_i = M$. Now let us define the model $M_{\omega}$. Let us choose some nontrivial ultrafilter $\mathbf{U}_{\omega}$ on $\omega$, and let the elements of the model $M_{\omega}$ be the sequences $(x_i)_{i<\omega}$ of elements of the model $M$, factorized by the ultrafilter $\mathbf{U}_{\omega}$. The model $M$ is embedded in $M_{\omega}$ by the constant sequences. Thus we have defined, in particular, some further extension of the field ${}^*\mathbf{R}$ of nonstandard real numbers. Now define $s_{\omega}$ as a class of sequences $(s_i)_{i<\omega}$ in $M_{\omega}$.

Every class $\lambda \in \Lambda$ is present in the enlarged model $M_{\omega}$, and in general the new classes appear. Let us denote the new set of classes by $\Lambda_{\omega} \supset \Lambda$. Choose a representative $a_{\lambda}$ from every class $\lambda \in \Lambda_{\omega}$, without changing the previously chosen ones. Now we may define $s_{\omega+1} = s_{\omega} + c_{\omega+1}a_{\lambda_{\omega+1}}$, etc. for all $i < 2\omega$. Then we pass to the new model $M_{2\omega}$; to do this, we choose a nontrivial
 ultrafilter $\mathbf{U}_{2\omega}$ on the set $\{i < 2\omega\}$ and construct the model $M_{2\omega}$ as the set of equivalence classes $(x_i)_{i<2\omega}$ of elements of the model $M_\omega$ modulo ultrafilter $\mathbf{U}_{2\omega}$. Then $s_{2\omega}$ is the class of the sequences $(s_i)_{i<2\omega}$, all terms of which are already constructed (the terms $s_i$ for $i < \omega$ belong to $M$, thus to $M_\omega$, because $M \subset M_\omega$). By the transfinite induction, we define the models $M_\beta$ for all $\beta < \alpha$; if $\beta = \gamma + 1$, then $M_\beta = M_\gamma$; if $\beta$ is a limit ordinal, then we introduce some nontrivial ultrafilter $\mathbf{U}_\beta$ on the ordinal $\beta$, and $M_\beta$ will consist of the equivalence classes of sequences $(x_i)_{i<\beta}$, $x_i \in M_i$, modulo ultrafilter $\mathbf{U}_\beta$. The term $s_\beta$ is defined as the class of the sequence $(s_i)_{i<\beta}$, if $\beta$ is a limit ordinal, and $s_\beta = s_\gamma + c_\beta a_{\lambda_\beta}$, if $\beta = \gamma + 1$.

The sets $\Lambda_\beta$ of equivalence classes are monotonely growing. Therefore, if we decompose the number $s_\alpha$ in the model $M_\alpha$, we obtain precisely the sequence $(c_i, \lambda_i)$. Thus, $M_\alpha$ is the required extension of the model $M$.

We may always take $\alpha$ sufficiently large ($|\alpha| > |\Lambda|$), taking $c_i = 0$ at the end of the sequence, so that we work with the same model $M_\alpha = M'$ for all sequences $(c_i, \lambda_i)$, $\lambda_i \in \Lambda$.

Now we may formulate our assertion.

\bigskip
\noindent\textbf{Theorem 4.2} \textit{The space $\mathcal{F}_M$ may be isometrically embedded into the asymptotic space $\mathcal{H}_{0,M'}$.}

\bigskip
\noindent\textbf{Proof.} For every sequence $(c_i, \lambda_i)_{i < \alpha < |\Lambda_M|}$, where $\lambda_i \in \Lambda_M$, there exists a number $s \in {{}^*\mathbf{R}_{M'}}$, such that $\operatorname{spec}(s) = (\lambda_i)_{i<\alpha}$, and $f_s(\lambda_i) = c_i$. Thus, for every $f \in \mathcal{F}_M$ we may define a point $(\rho, \varphi) \in {{}^*H_{\varepsilon,M'}}$, such that $f(\lambda)$ is defined for $\lambda_1 \leq \lambda \leq \lambda_0$, $\lambda_1 = \lambda(e^{(-\rho)})$ if and only if $\rho = -\log a_{\lambda_1}$ ($\rho \in {{}^*\mathbf{R}_M}$); and $\varphi$ is a (unique) number in ${}^*S^1_{M'}$, such that $\operatorname{spec}(\varphi) = (\lambda_i)_{i<\alpha}$, $f_\varphi(\lambda_i) = c_i$. This mapping is an isometric embedding, as we may see from the formulas (4.3) and (4.7).

\begin{flushright}
Q.E.D.
\end{flushright}

\section*{5 \quad Saturated models}

Our results are not complete: we have proved that for every model $M$ and every infinitesimal $\varepsilon \in {{}^*\mathbf{R}_M}$ there exist two spaces: The space $\mathcal{H}_{0,M}$ and the space $\mathcal{F}_M$. Our results are the following:

\noindent (1) $\mathcal{H}_{0,M}$ is isometrically embedded into $\mathcal{F}_M$;

\noindent (2) The model $M$ may be extended to a larger model $M'$, so that $\mathcal{F}_M$ is isometrically embedded into $\mathcal{H}_{0,M'}$.

The natural question is, whether there exists a model $M$, such that the spaces $\mathcal{H}_{0,M}$ and $\mathcal{F}_M$ \textit{coincide}. It turns out that such models exist, and they are the \textit{saturated models}.

Remind, what a saturated model is (see[6]). The first-order properties of a superstructure $\mathcal{S}$ are expressed by \textit{formulas}, written in a formal language $\mathcal{L}$ (see Section 1). A set $T$ of formulas with one free variable $x$, $T = \{R(x,c)\}$ (where $c$ is a parameter) is called \textit{finitely satisfiable}, if for every finite set of formulas $R(x,c_i)$ there exists an $x \in M$, such that $x$ satisfies all the formulas $R(x,c_i)$.

The model $M$ is called $\kappa$-saturated, where $\kappa$ is a cardinal, if for every finitely satisfiable set $T$ of formulas, such that $|T| < \kappa$, there exists $x \in M$, such that $x$ satisfies all formulas of $T$.

The model $M$ is called \textit{saturated}, if it is $\kappa$-saturated for $\kappa = |\overline{M}|$.

It is known, that there exist saturated models of ${}^*\mathbf{R}$ having every power $\aleph_i$, such that the Generalized Continuum Hypothesis holds in the power $\aleph_i$: $\aleph_{i+1} = 2^{\aleph_i}$.

Suppose $M$ is a saturated model of $\mathcal{S}$. Then the following is true:

\bigskip
\noindent\textbf{Theorem 5.1} \textit{For every sequence $(c_i, \lambda_i)_{i<\alpha}$, $c_i \in \mathbf{R}$, $\lambda_i \in \Lambda_M$, where $(\lambda_i)_{i<\alpha}$ is a decreasing sequence in $\Lambda_M$, $\alpha < |{{}^*\mathbf{R}}|$, there exists an $x \in {{}^*\mathbf{R}_M}$, such that $(\lambda_i)_{i<\alpha}$ is a spectrum of $x$, and $f_x(\lambda_i) = c_i$.}

\bigskip
\noindent\textbf{Proof.} (The idea of this proof has been proposed by Yoram Hirshfeld.) Let us construct the sequence $(s_i)_{i<\alpha}$. For $i < \omega$, set $s_i = c_0 a_{\lambda_0} + \cdots + c_i a_{\lambda_i}$.

\noindent Consider the set of formulas

\begin{equation*}
|x - s_i| < |a_{\lambda_i}| \tag{5.1}
\end{equation*}

\noindent for all $i < \omega$. This is a set of $|\omega| < |\underline{M}|$ finitely satisfiable formulas, so there exists an $x$, satisfying all of them.

Thus, the set $X_{(c_i,\lambda_i),i<\omega}$, defined in Section 3, is not empty, and we define

\begin{equation*}
s_\omega = a_{(c_i,\lambda_i),i<\omega}, \tag{5.2}
\end{equation*}

\noindent where $a_{(c_i,\lambda_i),i<\omega}$ is the minimal element of the set $X_{(c_i,\lambda_i),i<\omega}$ with respect to the order relation $\prec$, defined in Section 3.

Then we define

\begin{equation*}
s_{\omega+i} = s_\omega + \sum_{j=\beta,1}^{i} c_{\omega+j} a_{\lambda_{\omega+j}}. \tag{5.3}
\end{equation*}

The set of formulas $|x - s_i| < |a_{\lambda_i}|$, $i < 2\omega$ is finitely satisfiable, and thus $X_{(c_i,\lambda_i),i<2\omega} \neq \emptyset$; let us define

\begin{equation*}
s_{2\omega} = a_{(c_i,\lambda_i),i<2\omega}. \tag{5.4}
\end{equation*}

Proceeding by the transfinite induction, we define $s_i$ successively for all non-limit ordinals.

If $\beta < \alpha$, $\beta = \gamma + 1$, then we define $s_\beta = s_\gamma + c_\beta a_{\lambda_\beta}$.

If $\beta < \alpha$ is a limit ordinal, and $s_\gamma$ is defined for all $\gamma < \beta$, then the set of formulas

\begin{equation*}
|x - s_\gamma| < |a_{\lambda_\gamma}|, \quad \gamma < \beta, \tag{5.5}
\end{equation*}

\noindent is finitely satisfiable (for example, by $s_{\gamma_0}$, where $\gamma_0 = \max(\gamma_1, \cdots, \gamma_n)$, if we take a finite number of formulas, corresponding to $\gamma_1, \cdots, \gamma_n$). Thus, $X_{(c_i,\lambda_i),i<\beta}$ is not empty, and we define

\begin{equation*}
s_\beta = a_{(c_i,\lambda_i),i<\beta}. \tag{5.6}
\end{equation*}

Thus, we reach the index $\alpha$, because $|\alpha| < |\underline{M}|$, and define

\begin{equation*}
s_\alpha = a_{(c_i,\lambda_i),i<\alpha}, \tag{5.7}
\end{equation*}

for the sequence $(c_i, \lambda_i)$ finishes before $i$ reached the ordinal $\alpha$.

The procedure of constructing the spectrum of $s_\alpha$ shows immediately that

\begin{equation*}
\operatorname{spec}(s_\alpha) = (\lambda_i), \; i < \alpha, \tag{5.8}
\end{equation*}

\noindent and

\begin{equation*}
f_x(\lambda_i) = c_i. \tag{5.9}
\end{equation*}

\begin{flushright}
Q.E.D.
\end{flushright}

Thus, for any saturated model $M$, the correspondence between the points $(\rho, \varphi)$ of $\mathcal{H}_{0,M}$ and the functions $f \in \mathcal{F}_M$ is one-to-one, and $\mathcal{H}_{0,M}$ is isometric to $\mathcal{F}_M$ for arbitrary infinitesimal $\varepsilon \in {{}^*\mathbf{R}_M}$.

\section{The R-tree structure of the space $\mathcal{H}_{0,M}$ in the case of saturated model $M$}

The general principle is that for a hyperbolic space $\mathcal{X}$, its asymptotic space $\mathcal{X}_0$ is an $\mathbf{R}$-tree (see[5]). Let us look how this principle materialized in the case of the Lobachevsky plane $\mathcal{H}$.

\bigskip

\noindent Assume that the model $M$ is saturated.

\bigskip
\noindent\textbf{Theorem 6.1} $\mathcal{H}_0$ \textit{is a homogeneous metric space.}

\bigskip
\noindent\textbf{Proof.} We know that $\mathcal{H}$ is a homogeneous space; let $G$ be its transitive group of isometries.

\indent Using the transfer principle, we obtain that ${}^*\mathcal{H}$ is a homogeneous, with the isometry group ${}^*G$. (Isometries of ${}^*H$ are defined by the same explicit formulas).

\indent The space ${}^*\mathcal{H}_\varepsilon$ is also homogeneous, because we have applied the homothety to the homogeneous space. Denote by ${}^*G_\varepsilon$ the isometry group of ${}^*\mathcal{H}_\varepsilon$.

\indent Now observe that the isometries of ${}^*H_\varepsilon$, sending the origin to the finitely remote points of ${}^*\mathcal{H}_\varepsilon$, form a subgroup of ${}^*G_\varepsilon$; this is evident from the triangle inequality. So, $\mathcal{H}_0 = \mathrm{st}^*\mathcal{H}_\varepsilon$ is also homogeneous.

\begin{flushright}
Q.E.D.
\end{flushright}

\bigskip
\noindent\textbf{Theorem 6.2} \textit{If $M$ is a saturated model, then $\mathcal{H}_{0,M} = \mathcal{F}_M$ is a tree.}

\bigskip
\noindent\textbf{Proof.} We have to prove that the space $\mathcal{F}_M$ is connected, and every two elements may be connected by one single path (without self-intersections).

\indent Let $f_1, f_2$ be two elements of $\mathcal{F}_M$; these are two functions, defined on the segments $\lambda_{f_1} \leq \lambda \leq \lambda_0$, $\lambda_{f_2} \leq \lambda \leq \lambda_0$, respectively, of the ordered set $\Lambda$. These functions take their values in $S^1$ for $\lambda = \lambda_0$, and in $\mathbf{R}$ for other values of $\lambda$, and they are different from zero at the points $\lambda_{i,f_1}$ and $\lambda_{j,f_2}$, respectively, forming well-ordered sets with respect to the order on $\Lambda$ opposite to $<$. On $\Lambda$ there is defined a real-valued function $g(\lambda)$, distance from $\lambda_0$. This is a monotonically decreasing function, assuming all intermediate values (and all positive values). But every value $g$ of $g(\lambda)$ may be assumed on some segment $\Delta_g \subset \Lambda$. Let us choose a representative $\lambda_g \in \Delta_g$ for every $g \geq 0$.

\indent Define $h = \sup\{g | f_1(\lambda) = f_2(\lambda) \text{ for all } \lambda, \text{ such that} g(\lambda) \leq g\}$. Let $g_1 = g(\lambda_{f_1}),\ g_2 = g(\lambda_{f_2})$. Let $l = (g_1 - h) + (g_2 - h)$. For every $t \in [0, l]$, define

\begin{equation*}
g(t) = \begin{cases}
g_1 - t, & \text{if} \quad 0 \leq t \leq g_1 - h \\
g_2 - (l - t), & \text{if} \quad g_1 - h \leq t \leq l
\end{cases} \tag{6.1}
\end{equation*}

\noindent Now define the path $F_t$, $0 \leq t \leq l$ in $\mathcal{F}_M$:

\begin{equation*}
F_t = \begin{cases}
f_1|_{[\lambda_{g(t)},\lambda_0]}, & \text{if} \quad 0 \leq t \leq g - h \\
f_2|_{[\lambda_{g(t)},\lambda_0]}, & \text{if} \quad g_1 - h \leq t \leq l
\end{cases} \tag{6.2}
\end{equation*}

\noindent This is a continuous and isometric mapping of $[0, l]$ into $\mathcal{F}_M$. Thus, $\mathcal{F}_M$ is geodesically connected.

\indent It remains to prove that $\mathcal{F}_M$ is a tree. Suppose that $F_t'$, $0 \leq t \leq T$, is a path connecting $f_1$ and $f_2$. We prove that for some $t$, $F_t'$ is equal to $f_1$ restricted on $[\lambda', \lambda_0]$, such that $g(\lambda_0) = h$. Moreover, the path $F_t'$ should pass every point of the path $F_t$. It is evident, that this is enough. Suppose that the moments $t_i, 0 \leq i \leq N$ are so close that $\delta(F_{t_i}', F_{t_{i+1}}') < \delta_0$, where $\delta_0$ is small. Then the domains of $F_{t_i}'$ and $F_{t_{i+1}}'$ may differ not more than by $\delta_0$, and $F_{t_i}' \equiv F_{t_{i+1}}'$ outside the $\delta_0$-neighborhood of the endpoint of domain.

\indent If $g(\lambda_{t_i}) - h > c + \delta_0$ for all $i$, then the function $F_t'$ coincides with $f_1$ for $0 \leq g(\lambda) < h + c$; thus the distance $\delta(F_{t_i}', f_2) > g_2 - h + c$, and the path $F_t'$ cannot connect $f_1$ and $f_2$. Thus, for some $t$, $|g(\lambda_t) - h| < \delta_0$, and $F_t' \equiv f_1(t)$ for all $\lambda$, such that $g(\lambda) < h - \delta_0$. (We have to take $t = t_i$ with minimal $i$, such that $|g(\lambda_{t_i} - h| < \delta_0)$. This is true for all $\delta_0$; thus this is true for $\delta_0 = 0$, too. The same reasoning proves that $F_t'$ passes through every point of the path $F_s$ ($0 \leq s \leq g_1 - h$). By the symmetry, it is true also for all $s$, such that $g_1 - h \leq s \leq (g_1 - h) + (g_2 - h)$.

\begin{flushright}
Q.E.D.
\end{flushright}

\indent To conclude, note that the proof of homogeneity of the space $\mathcal{F}_M$ has nothing in common with the proof of homogeneity of the space $\mathcal{C}$, asymptotic subcone to $\mathcal{H}$, given in [8].

\bigskip
\section*{References}

\begin{enumerate}
\item[{[1]}] G. W. Brumfiel, \textit{The Tree of a Non-Archimedean Hyperbolic Plane}, in ``Geometry of Group Representations'' (Boulder, CO, 1987), Contemporary Mathematics, v. 74, 1988, p. 83-106.

\item[{[2]}] G. W. Brumfiel, \textit{The Real Spectrum Compactification of Teichm\"{u}ller Space,}, in ``Geometry of Group RepresentTions'' (Boulder, CO, 1987), Contemporary Mathematics, v. 74, 1988, p. 51-75.

\item[{[3]}] C. C. Chang and H. J. Keisler, \textit{Model Theory}, North Holland, Amsterdam-London-New York, 1973.

\item[{[4]}] M. Davis, \textit{Applied Nonstandard Analysis}, J. Wiley \& Sons, New York, 1977.

\item[{[5]}] M. Gromov, \textit{Asymptotic invariants of infinite groups// Geometric group theory. Vol. 2 (Sussex, 1991). London Mathematical Society Lecture Notes Ser. 182}, Cambridge University Press, Cambridge, 1993.

\item[{[6]}] M. Gromov \textit{Hyperbolic Groups// Essays in group theory/ ed. S. M. Gersten}, Springer, Berlin, 1987.

\item[{[7]}] A. H. Lightstone, A. Robinson \textit{nonarchimedean fields and asymptotic expansions}, North Holland, Amsterdam, 1975.

\item[{[8]}] J. Polterovich, A. Shnirelman, \textit{Asymptotic subcone of the Lobachevsky plane as a space of functions}, Uspehi Mat. Nauk, v. 52, number 4, p. 209-210, 1997.
\end{enumerate}

}\end{document}